# The Equation $\alpha \sin x + \beta \cos x = \gamma$ and a family of Heron Cyclic Quadrilaterals

*Konstantine Zelator*
*Department Of Mathematics*
*College Of Arts And Sciences*
*Mail Stop 942*
*University Of Toledo*
*Toledo,OH 43606-3390*
*U.S.A.*

# 1 Introduction

On the average, introductory texts in trigonometry may contain five to ten printed pages of material on trigonometric equations. Only in passing one may come across an equation such as $\sin x + \cos x = 1$; or, listed as an exercise in the exercise set. To study trigonometric equations and trigonometry in general, in more depth, one must resort to math books on advanced trigonometry. There are only a few of those around but they tent to be very good and thorough sources of information on the subject. Such are the two books listed in the references [1] and [2]. In [3] and [4], among other material, trigonometric equations (and systems of equations) are studied in some depth, and a variety of families of equations are analyzed. In [3], many types of problems of varying difficulty (both solved and unsolved) are presented. Let $\alpha, \beta, \gamma$ be fixed real numbers and consider the equation

$$\alpha \sin x + \beta \cos x = \gamma \qquad (1)$$

In Section 3, our analysis leads to the determination of the solution set S of Equation (1). In Section 3, via a simple construction process, we exhibit the angle $\theta$ that generates the solution set S of (1) in the case $\alpha^2 + \beta^2 = \gamma^2$ (see Figure 1) and with $\alpha, \beta, \gamma > 0$; the angle $\theta$ is part of an interesting quadrilateral $\Gamma B \Gamma_2 \Gamma_1$.

In Section 2 of this paper, we generate a family of quadrilaterals $\Gamma B \Gamma_2 \Gamma_1$ which have three integer side lengths, the fourth side length being rational, one integral diagonal length, on rational diagonal length, and with the four angles have rational tangent values and the four vertices lying on a circle; hence the term 'cyclic' in this paper's title. In Section 8, we will see that a certain subfamily of the above family; consists solely of Heron Quadrilaterals; quadrilaterals with integer side lengths and diagonals, as well as integral area.

## 2 The Solution Set of the Equation $\alpha \sin x + \beta \cos x = \gamma;\ \alpha, \beta, \gamma \in \mathrm{R}$

To find the solution set S (a subset of $R$), we use the well known half-angle formulas, which are valid for any angle or $x$ (typically measured in radians or degrees) not of the form $2k\pi + \pi; x \neq 2k\pi + \pi, k \in Z$ .

Note that as a simple calculation shows, if the real numbers of form $2k\pi + \pi$ are members of the solution set to (1), then it is necessary that $\beta + \gamma = 0$. So, under the restriction $x \neq 2k\pi + \pi$ , by using the above formulas and after a few algebra steps we obtain the equivalent equation

$$(\beta + \gamma)\tan^2\left(\frac{x}{2}\right) - 2\alpha \tan\left(\frac{x}{2}\right) + (\gamma - \beta) = 0 . \quad (2)$$

If $\beta + \gamma = 0$, then straight from (1) one sees that the reals of the form $2k\pi + \pi, k \in Z$ are solutions to (1). The other solutions are those reals which are solutions to (2). We have,

| |
|---|
| Suppose $\beta+\gamma=0$; that is ,$\gamma = -\beta$. Then the solution set $S$ of (1) is given by: |
| (i) $S = \mathrm{R}$, if $\alpha=\beta=\gamma=0$. |
| (ii) $S = \mathrm{S}_1$, if $\alpha=0$ and $\beta \neq 0$. |
| (iii) $S = \mathrm{S}_1 \bigcup S_2$, if $\alpha \neq 0$, |
| where |
| $S_1 = \{x \mid x \in \mathrm{R}, x = 2k\pi + \pi, k \in Z\}$ and |
| $S_2 = \{x \mid x \in \mathrm{R}, x = 2k\pi + 2\varphi, k \in Z\}$; |
| where $\varphi$ is the unique angle with $-\frac{\pi}{2} < \varphi < \frac{\pi}{2}$ and $\tan\varphi = -\frac{\beta}{\alpha}$. |

Now, assume $\beta+\gamma \neq 0$. Then, equation (2) is a quadratic equation in $\tan\left(\frac{x}{2}\right)$; and the trinomial $f(t)=(\beta+\gamma)t^2-2\alpha t+(\gamma-\beta)$ is a quadratic polynomial function with discriminant D given by

$$D=\left[(-2\alpha)\right]^2-4(\gamma+\beta)(\gamma-\beta)=4(\alpha^2+\beta^2-\gamma^2).$$

Clearly, when $D<0$, i.e., $\alpha^2+\beta^2<\gamma^2$ the trinomial $f(t)$ has two conjugate complex roots. When $D=0 \Leftrightarrow \alpha^2+\beta^2=\gamma^2$, the trinomial has the double-real root, $r=r_1=r_2=\frac{\alpha}{\beta+\gamma}$; and for $D>0$ (i.e., $\alpha^2+\beta^2>\gamma^2$), the trinomial $f(t)$ has two distinct real roots, namely We have (the reader is urged to verify),

> Suppose $\beta+\gamma \neq 0$; then the solution set $S$ of (1) is given by:
> (i) $S=\varnothing$, (empty set), if $\alpha^2+\beta^2<\gamma^2$.
> (ii) $S=\left\{x \mid x\in R, x=2k\pi+2\theta, k\in Z\right\}$, if $\alpha^2+\beta^2=\gamma^2$;
> where $\theta$ is a unique angle such that $-\frac{\pi}{2}<\theta<\frac{\pi}{2}$ and $\tan\theta=\frac{\alpha}{\beta+\gamma}$.
> (iii) $S=T_1\bigcup T_2$, if $\alpha^2+\beta^2>\gamma^2$;
> where $T_j=\left\{x \mid x\in R, x=2k\pi+2\theta_j\right\}$;
> j=1,2, and where $\theta_j$ is the unique angle such that $-\frac{\pi}{2}<\theta_j<\frac{\pi}{2}$
> and $\tan\theta_j=r_j$; with $\mathrm{r}_j=\frac{\alpha+(-1)^{j+1}.\sqrt{\alpha^2+\beta^2-\gamma^2}}{\beta+\gamma}$.

When $\alpha,\beta,\gamma$ are all positive there is an obvious geometric interpretation; $\alpha,\beta,\gamma$ can be the side lengths of a triangle, so the case $\alpha^2+\beta^2<\gamma^2$ corresponds to a triangle with the angle $\Gamma$ being obtuse. The case $\alpha^2+\beta=\gamma^2$ is the one of a right triangle with $\Gamma=90^\circ$. In the case of $\alpha^2+\beta>\gamma^2$, all three angles are acute.

# 3 The Case $\alpha^2+\beta^2=\gamma^2$ and a Geometric Construction

Let $\alpha,\beta,\gamma$ be positive real numbers such that $\alpha^2+\beta^2=\gamma^2$. Let us look at Figure 1. We start with the right triangle $\overset{\triangle}{\Gamma BA}$ with $|\overline{\Gamma B}|=\alpha, |\overline{\Gamma A}|=\beta, |\overline{AB}|=\gamma, \angle BA\Gamma=2\omega, \angle\Gamma BA=90^\circ-2\omega, \angle B\Gamma A=90^\circ$. The line segment $\overline{BA}$ lies on a straight line with divides the plane into two half-planes. In the half-plane opposite the point $\Gamma$, we draw at B the ray perpendicular to $\overline{BA}$ and we choose the point $\Gamma_2$ such that $|\overline{B\Gamma_2}|=|\overline{\Gamma B}|=\alpha$. We also extend the segment $\overline{BA}$ in the direction opposite from the point B and we pick the point $\Gamma_1$ such that $|\overline{A\Gamma}|=|\overline{A\Gamma_1}|=\beta$; thus, $|\overline{B\Gamma_1}|=\beta+\gamma$. From the isosceles triangle $\overset{\triangle}{\Gamma B\Gamma_2}$ we have $\angle B\Gamma\Gamma_2=\angle B\Gamma_2\Gamma=\varphi$. In the isosceles triangle, $\overset{\triangle}{\Gamma A\Gamma_2}$ it is clear $\angle A\Gamma\Gamma_1=\omega$. Furthermore, $\angle\Gamma B\Gamma_2=\angle\Gamma BA+90^\circ=90^\circ-2\omega+90^\circ=180^\circ-2\omega.$

Let $\angle B\Gamma_1\Gamma_2=\theta$. From the right triangle $\overset{\triangle}{B\Gamma_1\Gamma_2}$, it is clear that $\tan\theta=\dfrac{\alpha}{\beta+\gamma}$. (1), when $\alpha^2+\beta^2=\gamma^2$ (the condition $\beta+\gamma\neq 0$ is obviously satisfied since $\alpha,\beta,\gamma>0$). Now, in the isosceles triangle, $\angle\Gamma B\Gamma_2$, the sum of its three angles must equal $180^\circ: 2\varphi+(180^\circ-2\omega)=180^\circ \Leftrightarrow \varphi=\omega$, with proves that the four points $\Gamma, B, \Gamma_1, \Gamma_2$ lie on a circle with $\overline{\Gamma_1\Gamma_2}$ being a diameter, by virtue of $\angle\Gamma_2 B\Gamma_1=90^\circ$. Therefore, since these four points lie on a circle, we mush also have $\angle B\Gamma\Gamma_2 = \angle B\Gamma_1\Gamma_2$; that is, $\varphi = \theta$. Altogether,

$$\boxed{\varphi = \omega = \theta}$$

Note that angle $\angle\Gamma\Gamma_1\Gamma_2 = \omega + \theta = 2\omega$

And $\angle\Gamma B\Gamma_2 = 90^\circ - 2\omega + 90^\circ = 180^\circ - 2\omega$

And so, $\angle\Gamma\Gamma_1\Gamma_2 + \angle\Gamma B\Gamma_2 = 180^\circ$,

confirming that $\Gamma\Gamma_1\Gamma_2\mathrm{B}$ is a cyclic quadrilateral.

At this point, we make an observation, namely one that can also independently (from angle $\varphi$) show that $\omega = \theta$. Indeed, from the right triangle $\Gamma BA$, we have $0^\circ < 2\omega < 90^\circ; 0^\circ < \omega < 45^\circ$ and so $0 < \tan\omega < 1$. Also, $\tan 2\omega = \frac{\alpha}{\beta}$. Applying the double-angle identity $\tan 2\omega = \frac{2\tan\omega}{1-\tan^2\omega}$ leads to the equation $\alpha\tan^2\omega + 2\beta\tan\omega - \alpha = 0$. Using the quadratic formula combined with $\alpha^2 + \beta^2 = \gamma^2$, leads to the two possible solutions $\tan\omega = \frac{\gamma-\beta}{\alpha}, \frac{-(\beta+\gamma)}{\alpha}$; but only $\frac{\gamma-\beta}{\alpha}$ is an actual solution since $0 < \tan\omega < 1$. Hence, $\tan\omega = \frac{\gamma-\beta}{\alpha}$, but $\frac{\gamma-\beta}{\alpha} = \frac{\alpha}{\gamma+\beta}$ because $\alpha^2 + \beta^2 = \gamma^2$. Therefore, $\tan\omega = \frac{\alpha}{\beta+\gamma} = \tan\theta \Leftrightarrow (\text{under } 0<\theta,\omega<\frac{\pi}{2})\, \omega = \theta$.

Next, observe that since $\overline{\Gamma_2\Gamma_1}$ is a diameter, we must have $\angle\Gamma_2\Gamma\Gamma_1 = 90^\circ$ and $\angle\Gamma\Gamma_1\Gamma_2 = \omega + \theta = 2\omega$. If we set $x = \left|\overline{\Gamma\Gamma_1}\right|, y = \left|\overline{\Gamma\Gamma_2}\right|$, we have $x = \left|\overline{\Gamma_2\Gamma_1}\right|.\cos 2\omega$, and y= $\left|\overline{\Gamma_2\Gamma_1}\right|.\sin 2\omega$ and (from the triangle $\Gamma \overset{\triangle}{B}A$)

$\cos 2\omega = \dfrac{\beta}{\gamma}, \sin 2\omega = \dfrac{\alpha}{\gamma}$

We have,

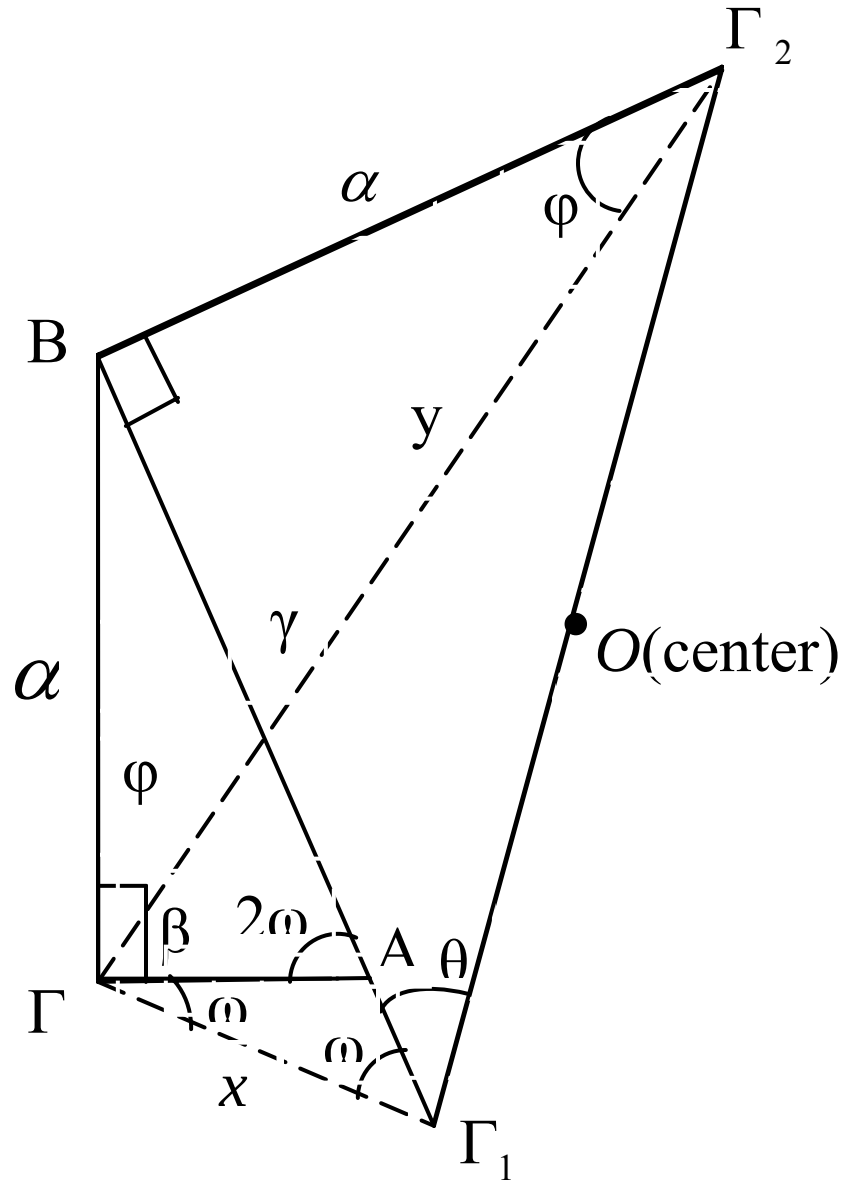

FIGURE 1

In the quadrilateral $\Gamma B\Gamma_2\Gamma_1$ in Figure 1,

(i) The four sidelengths are given by $\left|\overline{\Gamma B}\right| = \left|\overline{B\Gamma_2}\right| = \alpha$,

$\left|\overline{\Gamma_2\Gamma_1}\right| = \sqrt{\alpha^2 + (\beta+\gamma)^2}, \left|\overline{\Gamma\Gamma_1}\right| = x = \dfrac{\beta.\sqrt{\alpha^2 + (\beta+\gamma)^2}}{\gamma}.$

(ii) The two diagonal lengths are given by $\left|\overline{B\Gamma_1}\right| = \beta + \gamma$,

$\left|\overline{\Gamma\Gamma_2}\right| = y = \dfrac{\alpha.\sqrt{\alpha^2 + (\beta+\gamma)^2}}{\gamma}.$

(iii) The tangents of the four angles are given by, (3)

$$\tan(\angle\Gamma B\Gamma_2) = \tan(180^\circ - 2\omega) = -\tan 2\omega = -\frac{\alpha}{\beta},$$

$$\tan(\angle B\Gamma\Gamma_1) = \tan(90^\circ + \omega) = -\cot\omega = -\left(\frac{\alpha}{\gamma-\beta}\right) = \frac{\alpha}{\beta-\gamma},$$

$$\tan(\angle\Gamma\Gamma_1\Gamma_2) = \tan(\omega+\theta) = \tan 2\omega = \frac{\alpha}{\beta},$$

$$\tan(\angle\Gamma_1\Gamma_2 B) = \tan(90^\circ - \theta) = \cot\theta = \frac{\beta+\gamma}{\alpha}.$$

# 4 A numerical example

Let us take the perennial triple $(\alpha,\beta,\gamma) = (3,4,5)$. From (3) we obtain

$$\left|\overline{\Gamma_2\Gamma_1}\right| = \sqrt{90} = 3\sqrt{10}, \left|\overline{\Gamma\ \Gamma_1}\right| = x = \frac{4\sqrt{90}}{5} = \frac{12\sqrt{10}}{5}, \left|\overline{\Gamma B}\right| = \left|\overline{B\Gamma_2}\right| = 3, \left|\overline{B\Gamma_1}\right| = 9,$$

$$\left|\overline{\Gamma\Gamma_2}\right| = y = \frac{3\sqrt{90}}{5} = \frac{9\sqrt{10}}{5}, \tan\left(\angle\Gamma B\Gamma_2\right) = -\frac{3}{4}, \tan\left(\angle B\Gamma\Gamma_1\right) = -3, \tan\left(\angle\Gamma\Gamma_1\Gamma_2\right) = \frac{3}{4},$$

$$\tan\left(\angle\Gamma_1\Gamma_2 B\right) = 3.$$

From $\tan\theta = \dfrac{\alpha}{\beta+\gamma} = \dfrac{1}{3}$, with the aid of a scientific calculator we obtain

$\varphi = \omega = \theta \approx 18.43494882^\circ$ and by rounding up, $\varphi = \omega = \theta \approx 18.435^\circ$.

## 5 The case $\alpha, \beta, \gamma \in Z^+, \alpha^2 + \beta^2 = \gamma^2$ and a family of cyclic quadrilaterals with rational side lengths, rational diagonal lengths, and rational angle tangents.

When $\alpha, \beta, \gamma$ are positive integers such that $\alpha^2 + \beta^2 = \gamma^2$, then the triangle $\Gamma BA$ of Figure 1 is a Pythagorean triangle and $(\alpha, \beta, \gamma)$ is a Pythagorean triple. Note that in this case the real number $\sqrt{\alpha^2 + (\beta + \gamma)^2}$ will either be an irrational number or it will be a positive integer. It is an exercise in elementary theory to show if $n, c \in Z^+$ then $\sqrt[n]{c}$ will be a rational number if and only if, $c = k^n$, for some k $\in$ Z$^+$, so that $\sqrt[n]{c} = k$; we offer an explanation in Section 7 (see Fact 2), thus; the above square root will be rational, if and only if, it is a positive integer;

$$\sqrt{\alpha^2 + (\beta + \gamma)^2} = k \Leftrightarrow \alpha^2 + (\beta + \gamma)^2 = k^2, \text{ for some k} \in \text{Z}^+. \quad (4)$$

On the other hand, since $(\alpha, \beta, \gamma)$ is a Pythagorean triple (with $\gamma$ the hypotenuse length), we must have the following possibilities:

1. $\alpha = 2\delta mn, \beta = \delta(m^2 - n^2), \gamma = \delta(m^2 + n^2)$; or alternatively (5a)

2. $\alpha = \delta(m^2 - n^2), \beta = 2\delta mn, \gamma = \delta(m^2 + n^2)$ (5b)

for positive integers $\delta, m, n,$ such that m>n,(m,n)=1 (i.e., m and n are relatively prime) and m+n $\equiv$ 1(mod 2) (i.e., one of $m, n$ is odd, the other even). The above parametric formulas describe the entire family of Pythagorean triples. Derivation of these formulas can be found in most introductory number theory books (or texts); for example in [6].

We will only assume Possibility 1, i.e., (5a); and combine it with (3) to obtain a certain family of quadrilaterals $\Gamma B\Gamma_2\Gamma_1$; see note about Possibility 2 (i.e., (5b)) at the end of this section.

Combining (5a) and (4) yields,

$$4\delta^2 m^2(n^2+m^2)=k^2 \qquad (6)$$

It is also an exercise in elementary number theory to prove that if $a,b,n \in Z^+$ and $a^n$ is a divisor $b^n$, then $a$ must be a divisor of $b$. Again, refer to Section 7 of this paper, for an explanation (Fact 1). Thus, since by (6), the integer $4\delta^2 m^2=(2\delta m)^2$ is a divisor of $\mathrm{k}^2$. It follows that $k=2\delta mL,$ for some $L \in \mathrm{Z}^+$ and by (6),

$$\left.\begin{array}{l} m^2+n^2=L^2 \\ k=2\delta mL \end{array}\right\} \qquad (7)$$

We can now use (3), (5a), (6) and (7) to obtain expressions in terms of $m,n,\delta$ and $L$; of the four side lengths, the two diagonal lengths and the four tangent values of a special family of quadrilaterals:

Quadrilaterals $\Gamma B\Gamma_2\Gamma_1$

Family $F_1$:

(i) Sidelengths: $\left|\overline{\Gamma B}\right| = \left|\overline{B\Gamma_2}\right| = 2\delta mn,$

$\left|\overline{\Gamma_2\Gamma_1}\right| = 2\delta mL, x = \left|\overline{\Gamma\Gamma_1}\right| = \dfrac{2\delta m(m^2-n^2)}{L}$

(ii) Diagonal lengths: $\left|\overline{B\Gamma_1}\right| = 2\delta m^2, \left|\overline{\Gamma\Gamma_2}\right| = y = \dfrac{4\delta m^2 n}{L}$

$$\text{(iii) } \tan\left(\Gamma B\Gamma_2\right) = -\tan 2\omega = -\left(\frac{2mn}{m^2-n^2}\right) = \frac{2mn}{n^2-m^2} \qquad (8)$$

$$\tan\left(B\Gamma\ \Gamma_1\right) = -\cot\omega = -\left(\frac{2mn}{n^2}\right) = -\frac{2m}{n}$$

$$\tan\left(\Gamma\Gamma_1\Gamma_2\right) = \tan(\omega+\theta) = \tan 2\omega = \left(\frac{2mn}{m^2-n^2}\right)$$

$$\tan\left(\Gamma_1\Gamma_2 B\right) = \tan(90^\circ-\theta) = \cot\theta = \frac{2m}{n},$$

where $\delta, m, n, L$ are positive integers such that $m > n$, $(m,n) = 1, m+n \equiv 1(\text{mod } 2)$ and $m^2 + n^2 = L^2$.

Finally, in view of $(m,n)$=1 and the first equation in (7), we see that $(m,n,L)$ is itself a primitive Pythagorean triple which means,

$$\left.\begin{array}{l} m = t_1^2 - t_2^2, n = 2t_1t_2, L = t_1^2 + t_2^2 \text{ or alternatively} \\ m = 2t_1t_2, n = t_1^2 - t_2^2, L = t_1^2 + t_2^2 \end{array}\right\} \qquad (9a,9b)$$

where $t_1, t_2$ are positive integers such that $t_1 > t_2, (t_1, t_2) = 1, t_1 + t_2 \equiv 1(\text{mod}\,2)$; and always under the earlier assumption $m > n$.

**Note:** If one pursues Equation (5b) (Possibility 2) in combination with (6), one is led to the equation $\delta^2.(m+n)^2.2.(m^2+n^2) = k^2$ From there, using a little bit of elementary number theory on arrives at

$$\left.\begin{array}{l} m^2+n^2=2L^2 \\ k=2\delta(m+n)L \end{array}\right\}$$

This leads to a second family of quadrilaterals $\Gamma B\Gamma_2\Gamma_1$, which we will not consider here. We only point out that in this case, the triple $(m,n,L)$ will be a positive integer solution to the three variable Diophantine equation $X^2+Y^2=2Z^2$, whose general solution has been well known in the literature. The interested reader should refer to [5].

# 6 Another numerical example

If in (9b) we put $t_1=2, t_2=1$, then we obtain $L=5$ and m=4 > n=3 as required. Since the condition in (6) is satisfied if we take $\delta=5$, then we have, by (5a), $\alpha=120, \beta=35, \gamma=125$; and from (8) we obtain the quadrilateral $\Gamma B\Gamma_2\Gamma_1$ with the following specifications:

(i) Sidelengths: $\left|\overline{\Gamma B}\right|=\left|\overline{B\Gamma_2}\right|=120, \left|\overline{\Gamma_2\Gamma_1}\right|=200, x=\left|\overline{\Gamma\Gamma_1}\right|=56.$

(ii) Diagonal lengths: $\left|\overline{B\Gamma_1}\right|=160, y=\left|\overline{\Gamma\,\Gamma_2}\right| \;=92.$

(iii) $\tan\left(\Gamma B\Gamma_2\right)=-\frac{24}{7}, \tan\left(B\Gamma\Gamma_1\right)=-\frac{8}{3},$

$\tan\left(\Gamma\Gamma_1\Gamma_2\right)=\frac{24}{7}, \tan\left(\Gamma_1\Gamma_2 B\right)=\frac{8}{3}.$

Also, from $\tan\theta=\frac{\alpha}{\beta+\gamma}=\frac{24}{32}=\frac{3}{4}$, and with the aid of a scientific calculator, we find that $\varphi=\omega=\theta=36.86989765^\circ; \varphi=\omega=\theta=36.87^\circ$

**Remark:** Of the six lengths in (8), namely $\left|\overline{\Gamma B}\right|, \left|\overline{B\Gamma_2}\right|, \left|\overline{\Gamma_2\Gamma_1}\right|, x, \left|\overline{B\Gamma_1}\right|,$ and y, four are always integers as (8) clearly shows. These are the lengths $\left|\overline{\Gamma B}\right|, \left|\overline{B\Gamma_2}\right|, \left|\overline{\Gamma_2\Gamma_1}\right|$, and $\left|\overline{B\Gamma_1}\right|$.

On the other hand, the lengths $x$ and $y$ are rational but not integral unless $\delta$ is a multiple of $L$. This follows from the conditions $m+n \equiv 1 (\mathrm{mod}\, 2)$ and $(m,n)=1$. These two conditions and $m^2+n^2=L^2$ (see (7) or (8)), imply that the integer $L$ must be relatively prime or coprime to the product $2m(m^2-n^2)$ as well as to the product $4nm^2$; the proof of this is a standard exercise in an elementary number theory course. Therefore, according to (8), the rational number $x$ will be an integer precisely with $L$ is a divisor of $\delta$.

# 7 Two results from number theory

Let $a,b,n$ be positive integers.

**Fact 1:** If $a^n$ is a divisor of $b^n$, then a is a divisor of $b$.

**Fact 2:** The integer $a$ is the $n^{th}$ power of a rational number if and only if it is the $n^{th}$ power of an integer.

These two results can be typically found in number theory books. We cite two sources. First, W. Sierpinski's voluminous book *Elementary Theory of Numbers* (see reference [7] for details); and Kenneth H. Rosen's number theory text *Elementary Number Theory and Its Applications*, (see [6] for details).

# 8 A family of Heron Cyclic Quadrilaterals

If we compute the areas of the triangles $B\Gamma_1\Gamma_2, B\Gamma A$ and $\Gamma A\Gamma_1$ by using formulas (8) we find,

$$\left.\begin{aligned}
&\text{Area of (right) triangle } \overset{\square}{B\Gamma A} = \frac{\alpha\beta}{2} = \frac{2\delta mn\delta(m^2-n^2)}{2} = \delta^2.mn(m^2-n^2)\\
&\text{Area of (isosceles) triangle } \overset{\square}{\Gamma A\Gamma_1} = \frac{\beta^2.\sin(180^\circ-2\omega)}{2} = \frac{\beta^2\sin 2\omega}{2} = \frac{\beta^2}{2}\alpha/\gamma = \frac{\delta^2(m^2-n^2)^2 mn}{m^2+n^2}\\
&\text{Area of (right) triangle } \overset{\square}{B\Gamma_1\Gamma_2} = \frac{\alpha(\beta+\gamma)}{2} = \frac{2\delta mn\left[\delta(m^2+n^2)+\delta(m^2-n^2)\right]}{2} = 2\delta^2.n.m^3
\end{aligned}\right\} \quad (10)$$

The sum of the three areas in (10) is equal to the area A of the quadrilateral $B\Gamma\Gamma_1\Gamma_2$

$$A = \delta^2.m.n\left[m^2 - n^2 + \frac{(m^2-n^2)^2}{m^2+n^2} + 2m^2\right] \quad (11)$$

In the winter 2005 issue of Mathematics and Computer Education (see [8]), K.R.S. Sastry presented a family of Heron quadrilaterals. These are quadrilaterals with integer sides, integer diagonals, and integer area. Interestingly, a subfamily of the family of quadrilaterals described in (8); consists exclusively of heron quadrilaterals. First note that, for any choice of the positive integer $\delta$, the lengths $|\overline{\Gamma B}|, |\overline{B\Gamma_2}|, |\overline{\Gamma_2\Gamma_1}|, |\overline{B\Gamma_1}|$ are always integral, while $|\overline{\Gamma\,\Gamma_1}|$ and $|\overline{\Gamma\Gamma_2}|$, just rational. Clearly, if we take $\delta$ such that $\delta \equiv 0 \pmod L$, then $|\overline{\Gamma\,\Gamma_1}|, |\overline{\Gamma\Gamma_2}|$ will be integers as well; and as (11) easily shows, the area A will also be an integer, since by (7) $L^2 = m^2 + n^2$ and so $\delta \equiv 0 \pmod L$ => $\delta^2 \equiv 0 \pmod{L^2}$ => $\delta^2 \equiv 0 \pmod{m^2+n^2}$.

**Conclusion**: When $\delta$ is a positive integer multiple of $L$ in (8); the quadrilaterals obtained in (8) are Heron ones.

The smallest such choice for $\delta$ is $\delta$=L. From (8) we then obtain

$$\begin{aligned}
&|\overline{\Gamma\, B}| = |\overline{B\Gamma_2}| = 2Lmn, |\overline{\Gamma_2\Gamma_1}| = 2m.L^2, x = |\overline{\Gamma\,\Gamma_1}| = 2m(m^2-n^2),\\
&|\overline{B\Gamma_1}| = 2Lm^2, |\overline{\Gamma\,\Gamma_2}| = y = 4nm^2
\end{aligned} \quad (12)$$

And area $A$=L$^2$.$mn\left[ m^2 - n^2 + \frac{(m^2-n^2)^2}{m^2+n^2} + 2m^2 \right]$; and since $L^2 = m^2 + n^2$ we arrive at,

$$A = mn\left[ (m^2+n^2).\left(m^2-n^2\right) + (m^2-n^2)^2 + 2m^2.(m^2+n^2) \right] \quad (13)$$

$$A = mn\left[ m^4 - n^4 + m^4 - 2m^2n^2 + n^4 + 2m^4 + 2m^2n^2 \right] = 4nm^5$$

Using (12),(13), and (9a,9b) we arrive at the following table(with $t_1 t_2$ being the smallest possible; one choice with $t_1$ even $t_2$ odd; the other with $t_1$ odd, $t_2$ even).

| $t_1$ | $t_2$ | $m$ | $n$ | $\delta = L$ | $\left\|\overline{B\Gamma}\right\|$ | $\left\|\overline{\Gamma\Gamma_1}\right\|$ | $\left\|\overline{\Gamma_1\Gamma_2}\right\|$ | $\left\|\overline{\Gamma_2 B}\right\|$ | $\left\|\overline{B\Gamma_1}\right\|$ | $\left\|\overline{\Gamma\Gamma_2}\right\|$ | Area A |
|---|---|---|---|---|---|---|---|---|---|---|---|
| 2 | 1 | 4 | 3 | 5 | 120 | 56 | 200 | 120 | 160 | 192 | 12888 |
| 3 | 2 | 12 | 5 | 13 | 1560 | 2856 | 4056 | 1560 | 3744 | 2880 | $4.12^5.5$ |

Also note that from $\tan\theta = \frac{\alpha}{\beta+\gamma} = \frac{2\delta mn}{\delta(m^2-n^2)+\delta(m^2+n^2)} = \frac{n}{m}$;

we obtain (when $n$=3, $m$=4) $\tan\theta = \frac{3}{4}; \theta = \omega = \varphi \approx 36.86989765^\circ$;

and (when $n$=5, $m$=12) $\tan\theta = \frac{5}{12}; \theta = \omega = \varphi \approx 22.61986495^\circ$.